\documentclass[12pt]{article}
\usepackage[pctex32]{graphics}
\usepackage{amsthm,amsmath,amssymb,amscd,verbatim,epsfig}
\usepackage{amssymb}
\usepackage{graphicx}
\newcommand{\beq}{\begin{equation}}
\newcommand{\eeq}{\end{equation}}

\date{}

\newcommand{\sq}{$\square$}

\begin{document}

\title{More\ simple\ proofs\ of\ Sharkovsky's\ theorem}
\author{Bau-Sen Du \\ [.5cm]
Institute of Mathematics \\
Academia Sinica \\
Taipei 11529, Taiwan \\
dubs@math.sinica.edu.tw \\}
\maketitle

\section{Introduction}
Let $f$ be a continuous map from a compact interval $I$ into itself and let the Sharkovsky's ordering $\prec$ of the natural numbers be defined as follows:  $$3 \prec 5 \prec 7 \prec  9 \prec \cdots \prec 2 \cdot 3 \prec 2 \cdot 5 \prec 2 \cdot 7 \prec  2 \cdot 9 \prec \cdots \prec 2^2 \cdot 3 \prec 2^2 \cdot 5 \prec 2^2 \cdot 7 \prec  2^2 \cdot 9 \prec \cdots$$ $$\prec \cdots \prec 2^3 \prec 2^2 \prec 2 \prec 1.$$  Sharkovsky's theorem {\bf{\cite{mi}}}, {\bf{\cite{sh}}} states that (1) if $f$ has a period-$m$ point and if $m \prec n$ then $f$ also has a period-$n$ point; (2) for each positive integer $n$ there exists a continuous map $f : I \to I$ that has a period-$n$ point, but has no period-$m$ point for any $m$ with $m \prec n$; and (3) there exists a continuous map $f : I \to I$ that has a period-$2^i$ point for $i = 0, 1, 2, \ldots$ and has no periodic point of any other period.  There have been a number of different proofs {\bf [1-10, 12-17]} of it in the past 30 years providing various viewpoints of this beautiful theorem.  The sufficiency of Sharkovsky's theorem is well-known {\bf{\cite{du2}}} to be a consequence of the following three statements: (a) if $f$ has a period-$m$ point with $m \ge 3$, then $f$ also has a period-2 point; (b) if $f$ has a period-$m$ point with $m \ge 3$ and odd, then $f$ also has a period-$(m+2)$ point; and (c) if $f$ has a period-$m$ point with $m \ge 3$ and odd, then $f$ also has a period-$6$ point and a period-$(2m)$ point.  Note that in (c) we don't require the existence of periodic points of all even periods.  Only the existence of period-6 and period-$(2m)$ points suffices.  Among these three statements, (a) and (b) are easy to prove.  So, if we can find an easy proof of (c), then with the easy conterexamples in {\bf{\cite{al}}}, {\bf{\cite{du2}}} or {\bf{\cite{el}}}, we will have a simple proof of Sharkovsky's theorem.  In {\bf{\cite{du2}}}, {\bf{\cite{du3}}}, we present two different proofs of (c).  In this note, we present yet another four different proofs of (c) (including one that is a variant of the proof given in {\bf{\cite{du3}}}) and hence of Sharkovsky's theorem.

In the sequel, let $P$ be a period-$m$ orbit of $f$ with $m \ge 3$, let $b$ be the point in $P$ such that $f(b) = \min P$, and let $v$ be a point in $[\min P, b]$ such that $f(v) = b$.  

\section{A directed-graph proof of (c)}
This directed-graph proof of (c) is different from the one presented in {\bf{\cite{du2}}}.  Let $z_0 = \min \{ v \le x \le b : f^2(x) = x \}$.  Then $\min P = f^2(v) < v < z_0 < \max P$ and $f(x) > x$ for all $v \le x < z_0$.  Since $m \ge 3$ is odd, $\max P$ is also a period-$m$ point of $f^2$.  Hence $\max \{ f^2(x) : \min P \le x \le v \} > z_0$.  Let $I_0 = [\min P, v]$ and $I_1 = [v, z_0]$.  For each $n \ge 0$, the cycle $I_1(I_0)^nI_1$ (with respect to $f^2$) gives a period-$(n+1)$ point $w$ in $[v, z_0]$ for $f^2$ such that $f^{2i}(w) < w$ for all $1 \le i \le n$.  Since $f^{2i}(w) < w < f(w)$ for all $1 \le i \le n$, $w$ is a period-$(2n+2)$ point of $f$.  Therefore, $f$ has periodic points of all even periods.  This proves (c)

\section{A unified non-directed graph proof of (a), (b) and (c)}
The proof we present here is a variant of the proof given in {\bf{\cite{du3}}}.  Let $z$ be a fixed point of $f$ in $[v, b]$.  Then since $f^2(\min P) > \min P$ and $f^2(v) = \min P < v$, the point $y = \max \{ \min P \le x \le v : f^2(x) = x \}$ exists and $f(x) > z$ on $[y, v]$ and $f^2(x) < x$ on $(y, v]$.  So, $y$ is a period-2 point of $f$.  This proves (a).  On the other hand, assume that $m \ge 3$ is odd.  Then $f^{m+2}(y) = f(y) > y$ and $f^{m+2}(v) = \min P < v$.  So, the point $p_{m+2} = \min \{ y \le x \le v : f^{m+2}(x) = x \}$ exists and, since $f^2(x) < x$ on $(y, v]$, is a period-$(m+2)$ point of $f$.  This proves (b).  We now prove (c).  Let $z_0 = \min \{ v \le x \le b : f^2(x) = x \}$.  Then $z_0 \le z$ and $f(x) > z > x > f^2(x)$ on $[v, z_0)$.  If $f^2(x) < z_0$ on $[\min P, y]$, then $f^2(x) < z_0$ on $[\min P, z_0)$.  In particular, $f^2([\min P, z_0] \cap P) \subset [\min P, z_0] \cap P$.  Since $v$ lies in $[\min P, z_0]$ and since $f^2(v) = \min P$ is in $P$, we obtain that $f^{2n}(v) < z_0$ for all $n \ge 0$, contradicting the fact that $z_0 \le z < b = f^m(b) = f^m(f(v)) = f^{m+1}(v) = (f^2)^{(m+1)/2}(v)$.  Therefore, the point $d = \max \{ \min P \le x \le y : f^2(x) = z_0 \}$ exists and on $(d, y)$, we have $f(x) > z \ge z_0 > f^2(x)$.  Consequently, $f(x) > z \ge z_0 > f^2(x)$ on $(d, z_0)$.  For each $n \ge 1$, let $c_{2n} = \min \{ d \le x \le y : f^{2n}(x) = x \}$.  Then $d < \cdots < c_{2m+2} < c_{2m} < c_{2m-2} < \cdots < c_4 < c_2 \le y$.  Furthermore, $f^2(x) < z_0$ on $(d, c_2)$, $f^4(x) < z_0$ on $(d, c_4)$, $\cdots$, $f^{2n}(x) < z_0$ on $(d, c_{2n})$, $\cdots$, and so on.  In particular, for each $n \ge 1$, $f^{2i}(c_{2n}) < z_0 $ for all $0 \le i \le n-1$.  Since $f(x) > z \ge z_0 > f^2(x)$ on $(d, z_0)$, each $c_{2n}$ is a period-$(2n)$ point of $f$.  Therefore, $f$ has periodic points of all even periods.  (c) is proved.

\section{A prelimanary result}
In the previous section, when we prove (a), we obtain a side result.  That is, $f^2(v) < v < z < f(v) = b = f^m(b) = f^{m+1}(v)$.  So, when $m \ge 3$ is odd, we have $f^2(v) < v < z < (f^2)^{(m+1)/2}(v)$.  Surprisingly, these inequalities imply, by Lemma 1 below, the existence of periodic points of all periods for $f^2$.  However, the existence of periodic points of all periods for $f^2$ does not necessarily guarantee the existence of periodic points of all {\it even} periods for $f$.  It only guarantees the existence of periodic points of $f$ with least period $2k$ for each even $k \ge 2$ and least period $\ell$ or $2\ell$ for each odd $\ell \ge 1$.  We need a little more work to ensure the existence of periodic points of all even periods for $f$.  By doing some suitable "surgery" to the map $f$ to remove all those periodic points of $f$ of odd periods $j$ with $3 \le j \le m$, we can achieve our goal.  We present two such strategies in the next two sections.

\noindent
{\bf Lemma 1.}
If there exist a point $v$, a fixed point $z$ of $f$, and an integer $k \ge 2$ such that either $f(v) < v < z \le f^k(v)$ or $f^k(v) \le z < v < f(v)$, then $f$ has periodic points of all periods.  In particular, if $f$ has a periodic point of odd period $\ge 3$, then $f^2$ has periodic points of all periods.  

\noindent
{\it Proof.}
Without loss of generality, we may assume that $f(v) < v < z \le f^k(v)$ and $f$ has no fixed points in $[v, z)$.  Thus, $f(x) < x$ on $[v, z)$.  If $f(x) < z$ for all $\min I \le x \le v$, then $f(x) < z$ for all $x$ in $[\min I, z)$.  Consequently, $f^i(v) < z$ for each $i \ge 1$.  This contradicts the fact that $f^k(v) \ge z$.  So, the point $d = \max \{ \min I \le x \le v : f(x) = z \}$ exists and $f(x) < z$ on $(d, z)$.  Let $s = \min \{ d \le x \le z : f(x) \}$.  If $s > d$, then $d < s \le f(x) < z$ on $(d, z)$ which implies that $f^i(v) < z$ for each $i \ge 1$.  This, again, contradicts the fact that $f^k(v) \ge z$.  Thus, $s \le d$.  For each positive integer $n$, let $p_n = \min \{ d \le x \le z : f^n(x) = x \}$.  Then $p_n$ is a period-$n$ point of $f$.
\hfill\sq

\section{The second non-directed graph proof of (c)}
Let $z_1 \le z_2$ be the smallest and largest fixed points of $f$ in $[v, b]$ respectively.  Let $g$ be the continuous map on $I$ defined by $g(x) = \max \{ f(x), z_2 \}$ if $x \le z_1$; $g(x) = \min \{ f(x), z_1 \}$ if $x \ge z_2$; and $g(x) = -x + z_1 + z_2$ if $z_1 \le x \le z_2$.  Then $g([\min I, z_1]) \subset [z_2, \max I]$ and $g([z_2, \max I]) \subset [\min I, z_1]$ and $g^2(x) = x$ on $[z_1, z_2]$.  So, $g$ has no periodic points of any odd periods $\ge 3$.  Since $m \ge 3$ is odd, for some $1 \le i \le m-1$, both $f^i(b)$ and $f^{i+1}(b)$ lie on the same side of $[z_1, z_2]$.  Let $k$ be the smallest among these $i$'s.  If $k$ is odd then $g^2(v) = f^2(v) = f(b) = \min P < v < z_1 = g^{k+3}(v)$ and if $k$ is even then $g^2(v) = f^2(v) = f(b) = \min P < v < z_1 = g^{k+2}(v)$.  In either case, we have $g^2(v) < v < z_1 = (g^2)^n(v)$ for some $n \ge 2$.  By Lemma 1, $g^2$ has periodic points of all periods.  Since $g$ has no periodic points of any odd periods $\ge 3$, $g$ has periodic points of all even periods $\ge 4$ which obviously are also periodic points of $f$ with the same periods.  This proves (c).

\section{The third non-directed graph proof of (c)}
We may assume that $v$ is the largest point in $[\min P, b]$ such that $f(v) = b$.  Let $z_1 \le z_2$ be the smallest and largest fixed points of $f$ in $[v, b]$ respectively.  Without loss of generality, we may assume that $I = [\min P, \max P]$ and $f(x) = x$ for all $z_1 \le x \le z_2$.  By Lemma 1, we may also assume that $f$ has no fixed points in $[\min P, v] \cup [b, \max P]$.  If $f$ has a period-2 point $\hat y$ in $[v, z_1]$, then $z_2 < f(\hat y) < b$.  In this case, let $t = \hat y$ and $u = f(\hat y)$ and let $h$ be the continuous map from $I$ into itself defined by $h(x) = f(x)$ for $x$ not in $[t, u]$ and $h(x) = -x + t + u$ for $x$ in $[t, u]$.  Otherwise, let $t = v$ and $u = z_2$ and let $h(x) = f(x)$ for $x$ in $I$.  Consequently, if $Q$ is a period-$k$ orbit of $h$ with $k \ge 3$, then $[\min Q, \max Q] \supset [t, u]$.  Let $r = \min \{ \max Q : Q$ is a period-$m$ orbit of $h \}$.  Then $r$ is a period-$m$ point of $h$.  Let $R$ denote the orbit of $r$ under $h$.  By (b), $(\min R, \max R)$ contains a period-$(m+2)$ orbit $W$ of $h$.  Let $\hat h$ be the continuous map from $I$ into itself defined by $\hat h(x) = \min W$ if $h(x) \le \min W$; $\hat h(x) = \max W$ if $h(x) \ge \max W$; and $\hat h(x) = h(x)$ elsewhere.  Then clearly $\hat h$ has no period-$m$ points.  By (b), $\hat h$ has no period-$n$ points for any odd $n$ with $3 \le n \le m$.  Since $\hat h$ has the periodic orbit $W$ of odd period $m+2$, by Lemma 1 and the above, $\hat h$ has period-$(2j)$ points for each odd $j$ with $3 \le j \le m$ which obviously are also periodic points of $f$ with the same periods.  This proves (c).

\section{A proof of Sharkovsky's theorem}
For the sake of completeness, we include a proof of Sharkovsky's theorem which is slightly different from those in {\bf{\cite{du2}}}, {\bf{\cite{du3}}}.  If $f$ has period-$m$ points with $m \ge 3$ and odd, then by (b) $f$ has period-$(m+2)$ points and by (c) $f$ has period-$(2 \cdot 3)$ points.  If $f$ has period-$(2 \cdot m)$ points with $m \ge 3$ and odd, then $f^2$ has period-$m$ points.  According to (b), $f^2$ has period-$(m+2)$ points, which implies that $f$ has either period-$(m+2)$ points or period-$(2 \cdot (m+2))$ points.  If $f$ has period-$(m+2)$ points, then in view of (c) $f$ also has period-$(2 \cdot (m+2))$ points.  In either case, $f$ has period-$(2 \cdot (m+2))$ points.  On the other hand, since $f^2$ has period-$m$ points, by (c) $f^2$ has period-$(2 \cdot 3)$ points, hence $f$ has period-$(2^2 \cdot 3)$ points.  Now if $f$ has period-$(2^k \cdot m)$ points with $m \ge 3$ and odd and if $k \ge 2$, then $f^{2^{k-1}}$ has period-$(2 \cdot m)$ points.  It follows from what we have just proved that $f^{2^{k-1}}$ has period-$(2 \cdot (m+2))$ points and period-$(2^2 \cdot 3)$ points.  It follows that $f$ has period-$(2^k \cdot (m+2))$ points and period-$(2^{k+1} \cdot 3)$ points.  Consequently, if $f$ has period-$(2^i \cdot m)$ points with $m \ge 3$ and odd and if $i \ge 0$, then $f$ has period-$(2^\ell \cdot m)$ points for each $\ell \ge i$.  It is clear that $f^{2^\ell}$ has period-$m$ points.  By (a), $f^{2^\ell}$ has period-2 points.  This implies that $f$ has period-$2^{\ell+1}$ points whenever $\ell \ge i$.  Finally, if $f$ has period-$2^k$ points for some $k \ge 2$, then $f^{2^{k-2}}$ has period-4 points.  Again by (a), $f^{2^{k-2}}$ has period-2 points, ensuring that $f$ has period-$2^{k-1}$ points.  This proves the sufficiency of Sharkovsky's theorem.

For the rest of the proof, it suffices to assume that $I = [0, 1]$.  Let $T(x) = 1 - |2x - 1|$ be the tent map on $I$.  Then for each positive integer $n$ the equation $T^n(x) = x$ has exactly $2^n$ distinct solutions in $I$.  It follows that $T$ has finitely many period-$n$ orbits.  Among these period-$n$ orbits, let $P_n$ be the one with the smallest $\max P_n$ (or the one with the largest $\min P_n$).  For any $x$ in $I$ let $T_n(x) = \min P_n$ if $T(x) \le \min P_n$, $T_n(x) = \max P_n$ if $T(x) \ge \max P_n$, and $T_n(x) = T(x)$ elsewhere.  It is then easy to see that $T_n$ has exactly one period-$n$ orbit (i.e., $P_n$) but has no period-$m$ orbit for any $m$ with $m \prec n$ in the Sharkovsky ordering.  Now let $Q_3$ be the unique period-3 orbit of $T$ of smallest $\max Q_3$.  Then $[\min Q_3, \max Q_3]$ contains finitely many period-6 orbits of $T$.  If $Q_6$ is the one of smallest $\max Q_6$, then $[\min Q_6, \max Q_6]$ contains finitely many period-12 orbits of $T$.  We choose the one, say $Q_{12}$, of smallest $\max Q_{12}$ and continue the process inductively.  Let $q_0 = \sup \{\min Q_{2^i \cdot 3} : i \ge 0 \}$ and $q_1 = \inf \{ \max Q_{2^i \cdot 3} : i \ge 0 \}$.  Let $T_{\infty}(x) = q_0$ if $T(x) \le q_0$, $T_{\infty}(x) = q_1$ if $T(x) \ge q_1$, and $T_{\infty}(x) = T(x)$ elsewhere.  Then it is easy to check that $T_{\infty}$ has a period-$2^i$ point for $i = 0, 1, 2, \ldots$ but has no periodic point of any other period.  This completes the proof of Sharkovsky's theorem.

\end{document}